\documentclass[12pt,twoside,reqno,openany]{amsart}

\usepackage{amssymb,amscd,amsthm,amsmath,color,times}

\newcommand{\C}{\mathbb{C}}
\renewcommand{\P}{\mathbb{P}}
\newcommand{\X}{\mathbb{X}}

\author{Jo\"el Merker}

\title{Kobayashi hyperbolicity in degree $\geqslant n^{2n}$}

\begin{document}

\begin{abstract}
For a generic hypersurface $\mathbb{X}^{n-1} 
\subset \mathbb{P}^n(\mathbb{C})$ of degree
\[
d
\,\geqslant\,
n^{2n}
\]

\noindent
{\bf (1)}\,
$\mathbb{P}^n \big\backslash \mathbb{X}^{n-1}$ is 
Kobayashi-hyperbolically imbedded in $\mathbb{P}^n$;

\noindent
{\bf (2)}\,
$\mathbb{X}^{n-1}$ is 
Kobayashi($\Leftrightarrow$ Brody)-hyperbolic.

\smallskip

\noindent
{\bf (1)} improves Brotbek-Deng 1804.01719:
$d \geqslant (n+2)^{n+3}\, (n+1)^{n+3} = n^{2n}\,n^6\,
\big(e^3+{\rm O}(\frac{1}{n}) \big)$.

\noindent
{\bf (2)} supersedes Demailly 1801.04765:
$d \geqslant \frac{1}{3}\, \big( e^1(n-1) \big)^{2n}
= n^{2n}\, e^{2n}\, \big( \frac{1}{3\, e^2} + {\rm O}
(\frac{1}{n}) \big)$. 

\smallskip

\noindent
The method gives in fact 
$d \geqslant \frac{n^{2n}}{{\sf const}^n}$ for $n
\geqslant N({\sf const})$ with any 
${\sf const} \geqslant 1$.
\end{abstract}

\maketitle

\section{\bf Introduction}

The study of degeneracy properties of entire holomorphic curves
contained in projective algebraic varieties is a first (deep) step
towards understanding several outstanding problems of Diophantine
Geometry and Number Theory. In~\cite{Diverio-Merker-Rousseau-2010},
it was shown that such nonconstant entire curves in generic projective
hypersurfaces $\X^n \subset \P^{n+1}(\C)$ of degree $d \geqslant
2^{n^5}$ must land inside some fixed proper algebraic subvariety
$\mathbb{Y} \subsetneqq \mathbb{X}$, as predicted
by the Green-Griffiths conjecture. There, use was made of the
algebraic Morse inequalities of Trapani in order to
reduce the existence of global invariant jet differentials on
$\mathbb{X}$ to the positivity of a certain intersection number on the
Semple tower.  To handle the complexity and the difficulties of
computations not performed in~\cite{Diverio-Merker-Rousseau-2010}
within the cohomology ring of this tower, B\'erczi~\cite{Berczi-2015}
introduced equivariant localization to transform fixed point
formulas into iterated residue formulas which express various
intersection numbers between the tautological Semple line 
bundles, and he improved the degree bound to $d \geqslant n^{9n}$.

In his Ph.D. (Orsay, July 2014), Darondeau adapted and strengthened
these techniques to study strong algebraic degeneracy of entire curves
in complements $\P^n(\C) \backslash \X^{n-1}$ of generic projective
hypersurfaces, and in~\cite{Darondeau-2016} he reached $d \geqslant
(5n)^2\, n^n$. As is known by experts, Darondeau's calculations in the
complement case extend to the compact case, with the same degree
bound.

Siu~\cite{Siu-2015} developed a much more effective
method\,\,---\,\,concerning not only degree
bounds (also exponential), 
but mainly explicitness of jet differentials\,\,---\,\,to
establish the Kobayashi hyperbolicity conjecture.

Brotbek~\cite{Brotbek-2017} settled this conjecture
for $\mathbb{X}^n \subset \mathbb{P}^{n+1}(\C)$ 
(in a somehow
more general context), using a sufficiently rich subfamily of
hypersurfaces inspired from
Masuda-Noguchi \cite[7.4]{Noguchi-Winkelmann-2014} (after
explorations and exchanges with experts) on which he 
had the idea of using
Wronskian differential operators, associated multiplier ideal sheaves
and positivity of tautological line bundles on Grassmannian varieties.
From Brotbek's breakthrough, Deng~\cite{Deng-2016} deduced with short
arguments the degree bound $d \geqslant (n+1)^{n+2}\,
(n+2)^{2n+7} = n^{3n+9} \big( e^5 + {\rm O}(\frac{1}{n}) \big)$.

By means of an essentially smaller family of hypersurfaces
for which base loci of jet differentials can be
easily controlled, Demailly~\cite{Demailly-2018}
improved this to the lower degree
bound $d \geqslant \lfloor \frac{1}{3} \big(e\, 
n\big)^{2n+2} \rfloor$, still for generic
$\mathbb{X}^n \subset \mathbb{P}^{n+1}(\C)$.

Lowering the degree of examples of Kobayashi-hyperbolic
hypersurfaces is also a competitive topic.
Currently, the best existing examples in any dimension $n$
of such hyperbolic $\mathbb{X}^n \subset \mathbb{P}^{n+1}(\C)$
are due to Dinh-Tuan Huynh~\cite{Huynh-2016}:
they have any degree $d \geqslant \big(\frac{n+3}{2}\big)^2$. 
The reader is referred to this publication for more information.

Recently, by means of a general grassmannian technique having some
proximity with 
arguments of Shiffmann-Zaidenberg~\cite{Shiffman-Zaidenberg-2001}, 
Riedl-Yang~\cite{Riedl-Yang-2018} discovered that a solution
to the Green-Griffiths conjecture in dimension $2n$ implies a solution
to the Kobayashi conjecture in dimension $n$, with a small (linear)
degree bound discrepancy.  This motivates to put renewed efforts into
the Green-Griffiths conjecture for projective hypersurfaces.

\section{\bf Theorem and its proof}

Therefore, the goal of this concise note is to reveal that the
estimates of B\'erczi~\cite{Berczi-2015} and of Darondeau's Ph.D. 
({\em cf.}~\cite{Darondeau-2016}) can be slightly modified in order to
lower (improve) the current degree bounds for complex hyperbolicity.

\medskip\noindent{\bf Theorem.}
{\em For a generic hypersurface $\X^{n-1} \subset \P^n(\C)$ of degree}
\[
d
\,\geqslant\,
n^{2n}
\]

\smallskip\noindent{\bf (1)}\,
{\em $\P^n \big\backslash \X^{n-1}$ is Kobayashi-hyperbolically
imbedded in $\P^n$;}

\smallskip\noindent{\bf (2)}\,
{\em $\X^{n-1}$ is Kobayashi($\Leftrightarrow$ Brody)-hyperbolic.}

\medskip

\noindent
{\small\bf (1)} improves~\cite{Brotbek-Deng-2018}:
$d \geqslant (n+2)^{n+3}\, (n+1)^{n+3} = n^{2n}\,n^6\,
\big(e^3+{\rm O}(\frac{1}{n}) \big)$.

\noindent
{\small\bf (2)} supersedes~\cite{Demailly-2018}:
$d \geqslant \frac{1}{3}\, \big( e^1(n-1) \big)^{2n}
= n^{2n}\, e^{2n}\, \big( \frac{1}{3\, e^2} + {\rm O}
(\frac{1}{n}) \big)$. 

\smallskip

\noindent
The method gives in fact $d \geqslant \frac{n^{2n}}{{\sf const}^n}$
for $n \geqslant \text{\sc n}({\sf const})$ with any ${\sf const}
\geqslant 1$. We provide 
selected details only in the case ${\sf const} = 1$.

\proof
The way {\small\bf (1)} and {\small\bf (2)} are
related by induction on $n$ is known~\cite{Kobayashi-1998,
Noguchi-Winkelmann-2014}.
Treat only the complement case {\small\bf (1)}, since
Darondeau's modus operandi pp.~1894--1921 of~\cite{Darondeau-2016} 
applies mutatis
mutandis to the compact case. 
There, 
take $a_i := \big( \frac{n}{3} \big)^{n-i}$ for
$1 \leqslant i \leqslant n$\,\,---\,\,think
$\lceil \frac{n}{3} \rceil$; the choice
$a_i := \big( \frac{n}{\sf const} \big)^{n-i}$ would also work.
On p.~1915:
\[
\aligned
\widehat{C}\big(
{\textstyle{\frac{1}{a}}}
\big)
\,:=\,
&\,
\prod_{1\leqslant i<j\leqslant n}\,
\frac{a_i/a_j-1}{a_i/a_j-2}\,\,
\prod_{2\leqslant i<j\leqslant n}\,
\frac{a_i/a_j-2}{a_i/a_j-2-a_i/a_{i-1}}
\\
\,=\,
&\,
\prod_{k=1}^{n-1}\,
\bigg(
\frac{(n/3)^k-1}{(n/3)^k-2}
\bigg)^{n-k}\,\,
\prod_{k=1}^{n-1}\,\,
\bigg(
\frac{(n/3)^k-2}{(n/3)^k-2-\frac{1}{n/3}}
\bigg)^{n-1-k}
\\
\,=\,
&\,
\prod_{k=1}^{n-1}\,
\frac{\big((n/3)^k-1\big)^{n-k}}{
\big((n/3)^k-2\big)\,
\big(
(n/3)^k-2-\frac{1}{n/3}
\big)^{n-1-k}}.
\endaligned
\]
Introduce $C(\frac{1}{a})$ defined similary
by replacing $-\, a_i/a_{i-1}$ with $+a_i/a_{i-1}$:
\[
C\big(
{\textstyle{\frac{1}{a}}}
\big)
\,=\,
\prod_{k=1}^{n-1}\,
\frac{\big((n/3)^k-1\big)^{n-k}}{
\big((n/3)^k-2\big)\,
\big(
(n/3)^k-2+\frac{1}{n/3}
\big)^{n-1-k}}.
\]
The quotient:
\[
{\textstyle{\frac{\widehat{C}(\frac{1}{a})}{C(\frac{1}{a})}}}
\,=\,
\prod_{k=1}^{n-1}\,
\bigg(
\frac{
(n/3)^k-2+\frac{1}{n/3}
}{
(n/3)^k-2-\frac{1}{n/3}}
\bigg)^{n-1-k}
\,\,=\,\,\,
\prod_{k=1}^{n-1}\,
\left(
\frac{
1-\frac{2}{(n/3)^k}+\frac{1}{(n/3)^{k+1}}
}{
1-\frac{2}{(n/3)^k}-\frac{1}{(n/3)^{k+1}}}
\right)^{n-1-k}
\]
has logarithm asymptotic:
\[
\aligned
\log\,
{\textstyle{\frac{\widehat{C}(\frac{1}{a})}{C(\frac{1}{a})}}}
&
\,=\,
\sum_{k=1}^{n-1}\,
\big(n-1-k\big)\,
\bigg\{
\log\,
\Big[
1
-
\big(
{\textstyle{\frac{2}{(n/3)^k}}}
-
{\textstyle{\frac{1}{(n/3)^{k+1}}}}
\big)
\Big]
-
\log\,
\Big[
1
-
\big(
{\textstyle{\frac{2}{(n/3)^k}}}
+
{\textstyle{\frac{1}{(n/3)^{k+1}}}}
\big)
\Big]
\bigg\}
\\
&
\,=\,
\sum_{k=1}^{n-1}\,
\big(n-1-k\big)\,
\bigg\{
-
{\textstyle{\frac{1}{1}}}
\big(
{\textstyle{\frac{2}{(n/3)^k}}}
-
{\textstyle{\frac{1}{(n/3)^{k+1}}}}
\big)^1
-
{\textstyle{\frac{1}{2}}}
\big(
{\textstyle{\frac{2}{(n/3)^k}}}
-
{\textstyle{\frac{1}{(n/3)^{k+1}}}}
\big)^2
+
{\rm O}
\big(
{\textstyle{\frac{1}{n^3}}}
\big)
\,
+
\\
&
\ \ \ \ \ \ \ \ \ \ \ \ \ \ \ \ \ \ \ \ \ \ \ \ \ \ \ \
+
{\textstyle{\frac{1}{1}}}
\big(
{\textstyle{\frac{2}{(n/3)^k}}}
+
{\textstyle{\frac{1}{(n/3)^{k+1}}}}
\big)^1
+
{\textstyle{\frac{1}{2}}}
\big(
{\textstyle{\frac{2}{(n/3)^k}}}
+
{\textstyle{\frac{1}{(n/3)^{k+1}}}}
\big)^2
+
{\rm O}
\big(
{\textstyle{\frac{1}{n^3}}}
\big)
\bigg\}
\\
&
\,=\,
(n-2)\,
\Big\{
{\textstyle{\frac{1}{(n/3)^2}}}
+
{\textstyle{\frac{1}{(n/3)^2}}}
\Big\}
+
{\rm O}
\big(
{\textstyle{\frac{1}{n^2}}}
\big)
\,\,=\,\,
{\textstyle{\frac{18}{n}}}
+
{\rm O}
\big(
{\textstyle{\frac{1}{n^2}}}
\big).
\endaligned
\]
Similar estimations yield:
\[
\widehat{C}\big(
{\textstyle{\frac{1}{a}}}
\big)
\,=\,
e^3
+
{\rm O}
\big(
{\textstyle{\frac{1}{n}}}
\big).
\]

Next, define $I_0 := I_0^+ - I_0^-$ where as on pages 1896
and~1913:
\[
\aligned
I_0^\pm
&
\,:=\,
\pm\,
\sum_{k_1+\cdots+k_n\,=\,0
\atop
\pm\,C_{k_1,\dots,k_n}\,>\,0}\,
\big[
t_1^{n-k_1}\cdots t_n^{n-k_n}
\big]
\Big(
\big(
a_1t_1+\cdots+a_nt_n
\big)^{n^2}
\Big)
\cdot
C_{k_1,\dots,k_n},
\\
C^-
\big(
{\textstyle{\frac{1}{a}}}
\big)
&
\,:=\,
-\,
\sum_{k_1+\cdots+k_n\,=\,0
\atop
C_{k_1,\dots,k_n}\,<\,0}\,
C_{k_1,\dots,k_n}
\big(
{\textstyle{\frac{1}{a_1}}}
\big)^{k_1}
\cdots
\big(
{\textstyle{\frac{1}{a_n}}}
\big)^{k_n}
\,\,=\,\,
\frac{
\vert C\vert(\frac{1}{a})
-
C(\frac{1}{a})
}{2}.
\endaligned
\]
With $\widetilde{I}_0 := \frac{n^2!}{(n!)^n}\, a_1^n \cdots a_n^n$, 
on p.~1913:
\[
I_0^-
\,\leqslant\,
\widetilde{I}_0\,
C^-
\big(
{\textstyle{\frac{1}{a}}}
\big).
\]
On p.~1902:
\[
I_0^+
\,\geqslant\,
\widetilde{I}_0\,
\big(
1
+
3
+
{\textstyle{\frac{3^2}{2}}}
+
{\rm O}
({\textstyle{\frac{1}{n}}})
\big).
\]
From:
\[
\vert C\vert
\big(
{\textstyle{\frac{1}{a}}}
\big)
\big(
1
-
{\textstyle{\frac{18}{n}}}
-
{\rm O}({\textstyle{\frac{1}{n^2}}})
\big)
\,\leqslant\,
\widehat{C}
\big(
{\textstyle{\frac{1}{a}}}
\big)
\big(
1
-
{\textstyle{\frac{18}{n}}}
-
{\rm O}({\textstyle{\frac{1}{n^2}}})
\big)
\,\leqslant\,
C
\big(
{\textstyle{\frac{1}{a}}}
\big)
\,\leqslant\,
\vert C\vert
\big(
{\textstyle{\frac{1}{a}}}
\big)
\,\leqslant\,
\widehat{C}
\big(
{\textstyle{\frac{1}{a}}}
\big)
\,=\,
e^3
+
{\rm O}({\textstyle{\frac{1}{n}}}),
\]
it comes:
\[
\big(
0\,\leqslant\,
\big)
\ \ \ \ \ \ \ \ \ \ \ \ \
C^-
\big(
{\textstyle{\frac{1}{a}}}
\big)
\,\leqslant\,
{\textstyle{\frac{9\,e^3}{n}}}
+
{\rm O}({\textstyle{\frac{1}{n^2}}})
\,=\,
{\rm O}({\textstyle{\frac{1}{n}}}),
\]
hence:
\[
{\textstyle{\frac{1}{I_0}}}
\,\leqslant\,
{\textstyle{\frac{1}{\widetilde{I}_0}}}\,
\big(
{\textstyle{\frac{2}{17}}}
+
{\rm O}({\textstyle{\frac{1}{n}}})
\big).
\]
For $1 \leqslant p \leqslant n$ set as on p.~1898:
\[
\aligned
\widetilde{I}_p
\,:=\,
&\,
\widetilde{I}_0\,
\Big(
2n\,
\big(
1\,a_1
+\cdots+
n\,a_n
\big)
\Big)^p\,
\sum_{1\leqslant i_1<\cdots<i_p\leqslant n}\,
{\textstyle{\frac{1}{a_{i_1}\cdots\,a_{i_p}}}}
\\
\,=\,
&
\widetilde{I}_0\,
\Big[
{\textstyle{\frac{n^n}{3^n}}}\,
\big(
6
+
{\rm O}({\textstyle{\frac{1}{n}}})
\big)
\Big]^p\,
\Big(
\big(
{\textstyle{\frac{3}{n}}}
\big)^{\frac{p(p-1)}{2}}\,
\big(
1
+
{\rm O}({\textstyle{\frac{1}{n}}})
\big)
\Big).
\endaligned
\]
On p.~1906 using p.~1920:
\[
\aligned
\big\vert
I_p
\big\vert
&
\,\leqslant\,
\widetilde{I}_p
\cdot
\vert B\vert
\big(
{\textstyle{\frac{2n\mu h}{a_1}}},
\dots,
{\textstyle{\frac{2n\mu h}{a_n}}}
\big)
\cdot
\vert C\vert
\big(
{\textstyle{\frac{1}{a_1}}},
\dots,
{\textstyle{\frac{1}{a_n}}}
\big)
\cdot
{\textstyle{\frac{(n+1)^2+2}{2}}}
\\
&
\,\leqslant\,
\widetilde{I}_p
\cdot
\big(
e^{1/2}
+
{\rm O}({\textstyle{\frac{1}{n}}})
\big)
\cdot
\big(
e^3
+
{\rm O}({\textstyle{\frac{1}{n}}})
\big)
\cdot
{\textstyle{\frac{n^2}{2}}}
\big(
1
+
{\rm O}({\textstyle{\frac{1}{n}}})
\big),
\endaligned
\]
hence the largest root $\lambda(a)$ of $I_0 d^n + I_1 d^{n-1}
+ \cdots + I_n$ satisfies:
\[
\aligned
\lambda(a)
\,\leqslant\,
2\,
\underset{1\leqslant p\leqslant n}{\max}\,
\sqrt[p]{
{\textstyle{\frac{\vert I_p\vert}{I_0}}}}
\,\leqslant\,
2\,
\underset{1\leqslant p\leqslant n}{\max}\,
\Big(
{\textstyle{\frac{\widetilde{I}_p}{
\widetilde{I}_0}}}\,
{\textstyle{\frac{n^2}{2}}}\,e^{7/2}\,
{\textstyle{\frac{2}{17}}}
\big(
1
+
{\rm O}({\textstyle{\frac{1}{n}}})
\big)
\Big)^{\!\frac{1}{p}}
&
\,=\,
2\,
{\textstyle{\frac{\widetilde{I}_1}{\widetilde{I}_0}}}\,
n^2\,
{\textstyle{\frac{e^{7/2}}{17}}}
\big(
1
+
{\rm O}({\textstyle{\frac{1}{n}}})
\big)
\\
&
\,=\,
\frac{n^n}{3^n}\,
n^2\,
{\textstyle{\frac{12\,e^{7/2}}{17}}}\,
\big(
1
+
{\rm O}({\textstyle{\frac{1}{n}}})
\big).
\endaligned
\]
By~\cite{Riedl-Yang-2018}, {\small\bf (1)} holds for large
$n$ in degree:
\[
d
\,\geqslant\,
n^{2n}
\,\geqslant\,
n^{2n}\,
\big({\textstyle{\frac{2}{3}}}\big)^{2n}\,
(2n)^2\,
{\textstyle{\frac{12\,e^{7/2}}{17}}}\,
\big(
1
+
{\rm O}({\textstyle{\frac{1}{n}}})
\big),
\]
and for small $n$, a computer concludes, while
$d \geqslant n^n$ is not reachable so.
\endproof

\section{\bf Acknowledgments}

Thanks are adressed to Dinh-Tuan Huynh, Song-Yan Xie,
Mihai P\u{a}un for comments.


\begin{thebibliography}{10}

{\scriptsize

\bibitem{Berczi-2015}
B\'erczi, Gergely:
{\em Towards the Green-Griffiths-Lang conjecture via equivariant 
localisation},
{\tiny\sf arxiv.org/abs/1509.03406/}

\bibitem{Brotbek-2017}
Brotbek, Damian:
{\em On the hyperbolicity of general hypersurfaces}, 
Publ. Math. Inst. Hautes \'Etudes Sci. {\bf 126} (2017), 1--34. 

\bibitem{Brotbek-Deng-2018}
Brotbek, Damian; Deng Ya:
{\em Hyperbolicity of the complements 
of general hypersurfaces of high degree}, 
{\tiny\sf arxiv.org/abs/1804.01719/}

\bibitem{Darondeau-2016}
Darondeau, Lionel:
{\em On the logarithmic Green-Griffiths conjecture}, 
Int. Math. Res. Not. IMRN {\bf 2016}, no.~6, 1871--1923.
(Ph.D., Orsay, July 2014.)

\bibitem{Demailly-2018}
Demailly, Jean-Pierre:
{\em Kobayashi and Green-Griffiths-Lang conjectures}, 
{\tiny\sf arxiv.org/abs/1801.04765/}

\bibitem{Deng-2016}
Deng, Ya:
{\em Effectivity in the hyperbolicity-related problems}, 
{\tiny\sf arxiv.org/abs/1606.03831/}

\bibitem{Diverio-Merker-Rousseau-2010}
Diverio, Simone; Merker, Jo\"el; Rousseau, Erwan:
{\em Effective algebraic degeneracy},
Invent. Math. {\bf 180} (2010), no.~1, 161--223.

\bibitem{Huynh-2016}
Dinh~Tuan Huynh:
{\em Construction of hyperbolic hypersurfaces of low degree 
in $\mathbb{P}^n(\mathbb{C})$},
Internat. J. Math. {\bf 27} (2016), no.~8, 1650059, 9 pages.

\bibitem{Kobayashi-1998}
Kobayashi, Shoshichi:
Hyperbolic complex spaces. 
{\em Grundlehren der Mathematischen Wissenschaften}, 318. 
Springer-Verlag, Berlin, 1998. xiv+471~pp.

\bibitem{Noguchi-Winkelmann-2014}
Noguchi, Junjiro; Winkelmann, J\"org: 
Nevanlinna theory in several complex variables 
and Diophantine approximation. 
{\em Grundlehren der Mathematischen Wissenschaften}, 350. 
Springer, Tokyo, 2014. xiv+416~pp.

\bibitem{Riedl-Yang-2018}
Riedl, Eric; Yang, David: 
{\em Applications of a grassmannian technique in hypersurfaces}, 
{\tiny\sf arxiv.org/abs/1806.02364/}

\bibitem{Shiffman-Zaidenberg-2001}
Shiffman, Bernard; Zaidenberg, Mikhail:
{\em Hyperbolic hypersurfaces in $\mathbb{P}^n$ of Fermat-Waring type}, 
Proc. Amer. Math. Soc., {\bf 130} (2001), no.~7, 2031--2035.

\bibitem{Siu-2015}
Siu, Yum-Tong:
{\em Hyperbolicity of generic high-degree hypersurfaces 
in complex projective space}, 
Invent. Math. {\bf 202} (2015), no.~3, 1069--1166.

}

\end{thebibliography}
\end{document}